\documentclass[twoside]{article}
\usepackage{latexsym}
\usepackage{pstricks}
\usepackage{amssymb}
\usepackage{amsmath}
\usepackage{amsfonts}
\usepackage{amsthm,array}
\usepackage{array}
\usepackage{epsfig}
\usepackage{graphicx}

\def\bn{{\mathbb{N}}}

\def\br{{\mathbb{R}}}

\def\bz{{\mathbb{Z}}}

\def\br{\mathbb R}

\def\q{\quad}
\def\vs{\vskip.3cm}
\def\noi{\noindent}

\def\ve{\varepsilon}

\newrgbcolor{violet}{.6 .1 .8}
\newrgbcolor{lightyellow}{1 1 .8}
\newrgbcolor{lightblue}{.80 1 1}
\newrgbcolor{mygreen}{0 .66 .05}
\definecolor{mygreen}{rgb}{0,.66,.05}
\definecolor{lightyellow}{rgb}{1,1,.80}
\newrgbcolor{orange}{1 .6 0}
\newrgbcolor{GreenYellow}{.85 1 .31}
\newrgbcolor{Yellow}{1  1  0}
\newrgbcolor{Goldenrod}{1  .90  .16}
\newrgbcolor{Dandelion}{1  .71  .16}
\newrgbcolor{Apricot}{1  .68  .48}
\newrgbcolor{Peach}{1  .50  .30}
\newrgbcolor{Melon}{1  .54  .50}
\newrgbcolor{YellowOrange}{1  .58  0}
\newrgbcolor{Orange}{1  .39  .13}
\newrgbcolor{BurntOrange}{1  .49  0}
\newrgbcolor{Bittersweet}{1.  .4300  .24}
\newrgbcolor{RedOrange}{1  .23  .13}
\newrgbcolor{Mahogany}{1.  .4475  .4345}
\newrgbcolor{Maroon}{1.  .4084  .5376}
\newrgbcolor{BrickRed}{1.  .3592  .3232}
\newrgbcolor{Red}{1  0  0}
\newrgbcolor{OrangeRed}{1  0  .50}
\newrgbcolor{RubineRed}{1  0  .87}
\newrgbcolor{WildStrawberry}{1  .04  .61}
\newrgbcolor{Salmon}{1  .47  .62}
\newrgbcolor{CarnationPink}{1  .37  1}
\newrgbcolor{Magenta}{1  0  1}
\newrgbcolor{VioletRed}{1  .19  1}
\newrgbcolor{Rhodamine}{1  .18  1}
\newrgbcolor{Mulberry}{.6668  .1180  1.}
\newrgbcolor{RedViolet}{.9538  .4060  1.}
\newrgbcolor{Fuchsia}{.5676  .1628  1.}
\newrgbcolor{Lavender}{1  .52  1}
\newrgbcolor{Thistle}{.88  .41  1}
\newrgbcolor{Orchid}{.68  .36  1}
\newrgbcolor{DarkOrchid}{.60  .20  .80}
\newrgbcolor{Purple}{.55  .14  1}
\newrgbcolor{Plum}{.50  0  1}
\newrgbcolor{Violet}{.98 .15 .95}
\newrgbcolor{RoyalPurple}{.25  .10  1}
\newrgbcolor{BlueViolet}{.84  .38  .98}
\newrgbcolor{Periwinkle}{.43  .45  1}
\newrgbcolor{CadetBlue}{.38  .43  .77}
\newrgbcolor{CornflowerBlue}{.35  .87  1}
\newrgbcolor{MidnightBlue}{.4414  .9259  1.}
\newrgbcolor{NavyBlue}{.06  .46  1}
\newrgbcolor{RoyalBlue}{0  .50  1}
\newrgbcolor{Blue}{0  0  1}
\newrgbcolor{Cerulean}{.06  .89  1}
\newrgbcolor{Cyan}{0  1  1}
\newrgbcolor{ProcessBlue}{.04  1  1}
\newrgbcolor{SkyBlue}{.38  1  .88}
\newrgbcolor{Turquoise}{.15  1  .80}
\newrgbcolor{TealBlue}{.1572  1.  .6668}
\newrgbcolor{Aquamarine}{.18  1  .70}
\newrgbcolor{BlueGreen}{.15  1  .67}
\newrgbcolor{Emerald}{0  1  .50}
\newrgbcolor{JungleGreen}{.01  1  .48}
\newrgbcolor{SeaGreen}{.31  1  .50}
\newrgbcolor{Green}{0  1  0}
\newrgbcolor{ForestGreen}{.1992  1.  .2256}
\newrgbcolor{PineGreen}{.3100  1.  .5575}
\newrgbcolor{LimeGreen}{.50  1  0}
\newrgbcolor{YellowGreen}{.56  1  .26}
\newrgbcolor{SpringGreen}{.74  1  .24}
\newrgbcolor{OliveGreen}{.6160  1.  .4300}
\newrgbcolor{RawSienna}{.53  .28  .16}
\newrgbcolor{Sepia}{1.  .7510  .70}
\newrgbcolor{Brown}{.41  .25  .18}
\newrgbcolor{TAN}{.86  .58  .44}
\newrgbcolor{Gray}{1.  1.  1.}
\newrgbcolor{Black}{1  1  1}
\newrgbcolor{White}{1  1  1}
\newtheorem{theorem}{Theorem}[section]
\newtheorem{proposition}{Proposition}[section]
\newtheorem{lemma}{Lemma}[section]

\newtheorem{remark}{Remark}[section]

\newtheorem{remark-definition}{Remark and Definition}[section]
\newtheorem{rem-not}{Remark and Notation}[section]

\begin{document}
\title {Heteroclinic Orbits for a Discrete Pendulum Equation$^{\ast}$ }
\author {Huafeng Xiao  and  Jianshe Yu$^{\dagger}$ \\
{\small College of Mathematics and Information Sciences, Guangzhou University}\\
{\small Guangzhou 510405, People's Republic of China}}\date{}
\maketitle

 \footnote{\hskip.1cm$^{\ast}$ This project is supported by
National Natural Science Foundation of China (No. 10625104) and the
Research Fund for the Doctoral Program of Higher Education of China
(No. 20061078002).} \footnote{\hskip.1cm$^{\dagger}$ Corresponding
author. E-mail: jsyu@gzhu.edu.cn.}
\begin{abstract}
About twenty years ago, Rabinowitz showed firstly that there exist
heteroclinic orbits of autonomous Hamiltonian system joining two
equilibria. A special case of autonomous Hamiltonian system is the
classical pendulum equation. The phase plane analysis of pendulum
equation shows the existence of heteroclinic orbits joining two
equilibria, which coincide with the result of Rabinowitz. However,
the phase plane of discrete pendulum equation is similar to that of
the classical pendulum equation, which suggests the existence of
heteroclinic orbits for discrete pendulum equation also. By using
variational method and delicate analysis technique, we show that
there indeed exist heteroclinic orbits of discrete pendulum equation
joining every two adjacent points of $\{2k\pi+\pi:~k\in{\mathbb
Z}\}$. \vs

\end{abstract}

\section{Introduction}\label{sec:int}
Let us now introduce some notations that will be used throughout
this paper. By ${\mathbb {N,\,Z,\,R,\,R}^+}$ we denote the sets of
all natural numbers, integers, real numbers and positive real
numbers, respectively. For $a$, $b\in{\mathbb Z}$ ($a\le b$), define
{\it integer intervals} $Z[a]=\{a,a+1,a+2,\cdots\}$,
$Z[a,b]=\{a,a+1,\cdots,b\}$. For $D\subset{\mathbb R}$,
$\varepsilon>0$, denote by $B_{\varepsilon}(D)$ the open
$\varepsilon-$neighborhood of $D$. For a convergent bi-infinite
sequence $\{x_n\}_{n=-\infty}^{\infty}$, denote by $x_{\pm\infty}$
the limits of the sequence as $n$ tends to $\pm\infty$, i.e.
$x_{+\infty}:=\lim\limits_{n\rightarrow+\infty}x_n$ and
$x_{-\infty}:=\lim\limits_{n\rightarrow-\infty}x_n$.\vs

Consider the following second order equation
\begin{equation}\label{(1.1)}
\triangle^2x_{n-1}+A\sin x_n=0, \qquad\qquad  n\in {\mathbb Z},
\end{equation}
where $A\in{\mathbb R}^+$, $x_n\in{\mathbb R}$ for all $n\in{\mathbb
 Z}$, $\triangle$ is the forward difference operator
defined by $\triangle x_n=x_{n+1}-x_n$ and $\triangle^2
x_n=\triangle(\triangle x_n)$. A solution $x:{\mathbb Z}\rightarrow
{\mathbb R}$ of (\ref{(1.1)}) is called a {\it heteroclinic
solution} (or {\it heteroclinic orbit}) if there exist $\xi,
\eta\in{\mathbb R}$, $\xi\neq\eta$, such that $\xi,~\eta$ are two
equilibria of (\ref{(1.1)}) and $ x_{-\infty}=\xi,~
x_{+\infty}=\eta. $ \vs

We are interested in  the problem of the existence and multiplicity
heteroclinic solutions of (\ref{(1.1)}). So far as we are aware, it
is the first time in the literature for us to study heteroclinic
orbits of difference equations. \vs

Equation (\ref{(1.1)}) can be considered as a discrete analogue of
the classical pendulum equation:
\begin{equation}\label{(1.2)}
x^{\prime\prime}(t)+A\sin x(t)=0,\qquad t\in{\mathbb R}.
\end{equation}
The phase plane portrait of  (\ref{(1.1)}) with $A=0.1$, shown on
Figure \ref{fig:DE}, can be compared to the phase plane portrait of
(\ref{(1.2)}) with $A=1$, shown on Figure \ref{fig:ODE}. The phase
plane analysis of (\ref{(1.2)}) shows the existence of two
heteroclinic solutions for (\ref{(1.2)}) joining $-\pi$ and $\pi$.
On the other hand, the phase plane of (\ref{(1.1)}) is similar to
that of (\ref{(1.2)}). On Figure \ref{fig:DE}, we use  colors to
distinguish between different orbits. Nine ellipses represent nine
periodic orbits, while two curves around nine ellipses are
non-periodic orbits. Close similarities observed on Figures
\ref{fig:DE} and \ref{fig:ODE} suggest the existence of heteroclinic
orbits for (\ref{(1.1)}). Our goal in this paper is to show that
there indeed exist two heteroclinic solutions of (\ref{(1.1)})
joining $-\pi$ and $\pi$ also. \vs

Let us now recall briefly the existence and multiplicity
heteroclinic orbits for the following Hamiltonian system
\begin{equation}\label{(1.3)}
q^{\prime\prime}+V_q^{\prime}(t,q)=f(t),
\end{equation}
which is a generalization form of (\ref{(1.2)}), where
$q=(q_1,q_2,\dots,q_n)\in{\mathbb R}^n$, $V:{\mathbb
R}\times{\mathbb R}^n\rightarrow{\mathbb R}$ and $ f:{\mathbb
R}\rightarrow{\mathbb R}^n$. In the past twenty years, many authors
had studied the existence and multiplicity of heteroclinic solutions
and heteroclinic chains for (\ref{(1.3)}). The first result in this
area was proved in \cite{rabinowitz4}, where the author discussed
(\ref{(1.3)}) under assumptions: $(V_0)$ $f(t)=0$, $(V_1)$ $V\in
C^1({\mathbb R}^n,{\mathbb R})$, and $(V_2)$ $V$ is periodic in
$q_i$ with the period $T_i$, $1\leq i\leq n$. Conditions $(V_1)$ and
$(V_2)$ imply that $V$ has a global maximum on ${\mathbb R}^n$.
Without loss of generality, it is assumed that the global maximum of
$V$ is 0 and put $\mathbf{\Gamma}=\{\xi\in{\mathbb R}^n:\
V(\xi)=0\}$. Under the non-degeneracy condition: $(V_3)$
$\mathbf{\Gamma}$ consists only of isolated points, the following
result was obtained in \cite{rabinowitz4}. \vs

\noindent \textbf{Theorem A.} {\it Under assumptions $(V_0)-(V_3)$,
for every $\xi\in\mathbf{\Gamma}$, there exist at least two
heteroclinic orbits of (\ref{(1.3)}) joining $\xi$ to
$\mathbf{\Gamma}\setminus\{\xi\}$. At least one of these orbits
emanates from $\xi$ and at least one terminates at $\xi$}. \vs

Let $n=1$, $V(t,x)=A\sin x$ and $f(t)=0$. Then (\ref{(1.3)}) becomes
(\ref{(1.2)}), and $\mathbf{\Gamma}=\{2k\pi+\pi:~k\in{\mathbb Z}\}$.
Theorem A guarantees at least two heteroclinic orbits of
(\ref{(1.2)}) through every point of $\mathbf{\Gamma}$. \vs

Further development in this direction was done by Felmer (cf.
\cite{felmer}), who generalized the above result to first order
spatially periodic Hamiltonian systems. By using the saddle point
theorem, the author obtained the existence of heteroclinic orbits
joining two saddle type critical points. Without imposing the
non-degeneracy condition, Caldiroli and Jeanjean (cf.
\cite{caldiroli}) studied conservative singular Hamiltonian system
without forcing term, and were able to establish the existence of
heteroclinic orbits joining a global maximum point and a
non-constant periodic solution. \vs

In the case when the potential is periodic and time reversible, by
using minimization arguments, Rabinowitz showed the existence of
heteroclinic solutions between pairs of periodic solutions, (cf.
\cite{rabinowitz1,rabinowitz2}). Under the same assumptions,
 Maxwell (cf. \cite{maxwell}) proved that there exist heteroclinic chains connecting every pairs of
periodic solutions. \vs

For non-autonomous Hamiltonian systems without forcing term, Strobel
(cf. \cite{strobel}) studied the existence of heteroclinic chains
between pairs of equilibria. By using constrained minimization and
comparison arguments, Rabinowitz and Zelati (cf. \cite{rabinowitz5})
studied (\ref{(1.3)}) without forcing term, and found multiple
heteroclinic chains joining pairs of equilibria. Subsequent progress
was done by Bertotti and Montecchiari (cf. \cite{bertotti}), who
generalized the results of Strobel, by proving the existence of
infinitely many heteroclinic solutions for a class of forced slowly
oscillating Hamiltonian system with potential $V(t,x)$ of form
$\alpha(t)W(x)$, with $\alpha$ being almost periodic in $t$. Next,
Alessio, Bertotti and Montecchiari (cf. \cite{alessio}) obtained a
generalization of the results of \cite{bertotti}, in which
$\alpha(t)$ is replaced by $\alpha(t)+\alpha(\varepsilon t)$ for
$\varepsilon>0$ small enough. However, without the non-degeneracy,
these results are not as strong as those of Strobel. In the case of
forced slowly oscillating Hamiltonian systems, Rabinowitz (cf.
\cite{rabinowitz3}) showed the existence of basic and even more
complex heteroclinic orbits, without making any non-degeneracy
assumption. Then, Zelati and Rabinowitz (cf. \cite{zelati}) showed
that there exist heteroclinic solutions joining two stationary
points in different energy levels, which was established by using
minimization arguments. \vs

We should also mention the work by Chen and Tzend (cf.
\cite{chen1,chen2,chen3}), in which variational and penalization
methods were being used to study the existence of heteroclinic
orbits for the following system
\begin{equation}\label{(1.4)}
q^{\prime\prime}-V_q(t,q)=0,
\end{equation}
where $V$ is not periodic nor asymptotically periodic in $t$. In
those papers, the authors obtained multiple heteroclinic orbits and
chains joining pairs of equilibria as well as joining an equilibrium
to a non-constant periodic solution. Izydorek and Janczewska (cf.
\cite{izydorek}) proved, without assuming periodicity or almost
periodicity in $t$ of the potential, for (\ref{(1.3)}) without
forcing term, the existence of heteroclinic solutions joining pairs
of equilibria. \vs

However, no results on the existence of heteroclinic solutions of
difference equations have been proved. In this paper, by using
variational arguments, we will study the existence and multiplicity
of heteroclinic solutions for (\ref{(1.1)}). To this end, we need to
choose a suitable functional space on which a variational
functional, associated with (\ref{(1.1)}), can be constructed.
However, lack of continuity assumption leads to some new problems
which were not present in the case of differential systems. For
example, for differential systems, if an orbit contains two points
such that one of them is outside of $B_{\varepsilon}(\xi)$, while
the other belongs to inside of
$B_{\delta}(\xi)$($\delta<\varepsilon/2$), then the orbit (because
of its continuity) contains a point belonging to $\partial
B_{\varepsilon/2}(\xi)$. However, such a statement is not valid for
orbits of discrete systems. \vs

\section{Main Results}
In this section, we study the existence and multiplicity of
heteroclinic orbits of (\ref{(1.1)}) joining every two adjacent
points of $\{2k\pi+\pi:~k\in{\mathbb Z}\}$. For simplicity, we make
an image translation. By applying the substitution $y_n=x_n-\pi$,
(\ref{(1.1)}) can be rewritten as
\begin{equation}\label{(2.1)}
\triangle^2y_{n-1}-A\sin y_n=0, \qquad\qquad  n\in {\mathbb Z} .
\end{equation}
We look for heteroclinic orbits of (\ref{(2.1)}) which join two
adjacent points of $\{2k\pi:k\in{\mathbb Z}\}$. \vs

Let $C$ be the vector space of all convergent sequences
$y=\{y_k\}_{k=-\infty}^{\infty}$, i.e.
\begin{equation*}
C:=\Big\{y=\{y_k\}: \lim\limits_{k\rightarrow\infty}y_k  \
\mbox{and}\
   \lim\limits_{k\rightarrow-\infty}y_k\ \mbox{exist},
   ~y_k\in{\mathbb R}, ~k\in{\mathbb Z}\Big\}.
\end{equation*}
We define the space $H$ by
\begin{equation*}
H:=\left\{y\in C: \sum\limits_{k=-\infty}^{\infty}|\triangle
y_k|^2<\infty\right\},
\end{equation*}
and put
\begin{align}
<x,y>&:=\sum\limits_{k=-\infty}^{\infty}\triangle x_k\triangle
y_k+x_0y_0,\qquad \forall\ x,y\in H,\label{inner}\\
\|y\|&:=\left[\sum\limits_{k=-\infty}^{\infty}(\triangle
y_k)^2+y_0^2\ \right]^\frac{1}{2},\qquad \forall\ y\in
H.\label{norm}
\end{align}
Then we have \vs
\begin{proposition}\label{1}
The bilinear product (\ref{inner}) is an inner product on $H$ and
the space $H$ equipped with the norm given by (\ref{norm}) is a
Hilbert space.
\end{proposition}
\begin{proof} Recall that the space $l^2(\bz)$ of all sequences
 $a=\{a_k\}_{k=-\infty}^{\infty}$ such that
\[ \|a\|_2:=\left[\sum_{-\infty}^\infty a_k^2\right]^{\frac 12}<\infty,
\]
is a Hilbert space. Let $\{y^n\}\subset H$ be a Cauchy sequence in
$H$, i.e.
\begin{equation}\label{Cauchy}
\forall_{\ve>0}\; \exists_{N}\;\forall_{m,n\ge N}\;\;\;
 \|y^n-y^m\|=\left[\sum_{k=-\infty}^\infty (\triangle y_k^n-\triangle
 y_k^m)^2+(y_0^n-y_0^m)^2   \right]^{\frac 12}<\ve.
\end{equation}
Then $\{y_0^n\}$ is a Cauchy sequence in $\br$, while $\{\triangle
 y^n\}$,  $\triangle y^n:=\{\triangle y_k^n\}$, is a Cauchy sequence in
 $l^2(\bz)$. By completeness of $l^2(\bz)$, there exists a limit $a$ in
 $l^2(\bz)$ of $\{\triangle y^n\}$. One can easily observe, that there
 exists a unique $y^0:=\{y^0_k\}$ in $H$ such that
\[\lim_{n\to \infty} y^n_0=y^0_0,\quad \text { and }\q \forall_{k\in
 \bz} \;\; \triangle y^0_k=a_k.\]
By passing to the limit as $m$ goes to $\infty$, we obtain from
 (\ref{Cauchy})
 \[
 \forall_{\ve>0}\; \exists_{N}\;\forall_{n\ge N}\;\;\;
 \|y^n-y^0\|\le\ve,
 \]
 which proves that $\{y^n\}$ converges to $y^0$. Consequently,  $H$ is
 a Hilbert space.
\end{proof}
\vs

Similar arguments as those presented in  \cite{guo}, one can define
variational functional $J:H\to \overline{\br}:=\br\cup\{\infty\}$
associated with (\ref{(2.1)}) by
\begin{equation}\label{(2.3)}
J(y)=\sum\limits_{n=-\infty}^{\infty}[\frac{1}{2}|\triangle
y_n|^2+A(1-\cos y_n)].
\end{equation}
Put $\mathbf{\Theta}:=\{2k\pi:~k\in {\mathbb Z}\}$ and
$\gamma:=\frac{2\pi}{3}$.
\begin{remark}\label{R3}\rm
For every $y=\{y_n\}\in H$, if $J(y)<\infty$, then $y_{-\infty}$,
$y_{+\infty}\in\mathbf{\Theta}$. Indeed, suppose for example
$y_{+\infty}\notin \mathbf \Theta$, then  there exist
$\gamma>\delta>0$ and $N\in \mathbb N$ such that $y_n\notin
B_{\delta}(\mathbf{\Theta})$ for all $n\ge N$. Therefore,
\begin{eqnarray}
J(y)\geq\sum\limits^{\infty}_{i=N}A(1-\cos
y_i)\geq\sum\limits^{\infty}_{i=N}A(1-\cos \delta)=\infty.\nonumber
\end{eqnarray}
\end{remark}

Given $\xi\in \mathbf{\Theta}\setminus\{0\}$,
$\varepsilon\in(0,\gamma)$, define the set
$\Gamma_{\varepsilon}(\xi)$ of all $y\in H$ satisfying
\begin{itemize}
\item[(i)]  $ y_{-\infty}=0$,
\item[(ii)]  $y_{+\infty}=\xi$,
\item[(iii)] $y_n\notin B_{\varepsilon}(\mathbf{\Theta}\setminus\{0,\xi\})$ for
 all $n\in{\mathbb Z}$.
\end{itemize}
Obviously, $\Gamma_{\varepsilon}(\xi)$ is not empty for all
$\xi\in\mathbf{\Theta}$. Define
\begin{equation}\label{c-alpha}
c_{\varepsilon}(\xi):=\inf\limits_{y\in\Gamma_{\varepsilon}(\xi)}
J(y)\quad \mbox{and}\quad
 \alpha_{\varepsilon}:=\min\limits_{t\notin
 B_{\varepsilon}(\mathbf{\Theta})}(1-\cos t)>0.\nonumber
\end{equation}
Now we give a simple but useful lemma.
\begin{lemma}\label{2}
Given a sequence of disjoint integer intervals $Z(n_n,m_k)$,
$n_k<m_k$
 and $j\in \mathbb N$. Let $y\in H$ be such that
\[ y_i\notin B_{\varepsilon}(\mathbf{\Theta})\q \text{for}\q i\in
\bigcup_{k=1}^j Z(n_k,m_k).\]
 Then,
\[
J(y)\geq\sqrt{{2A\alpha_{\varepsilon}}}\sum\limits_{k=1}^j|y_{m_k}-y_{n_k}|.
\]
\end{lemma}
\begin{proof}
Let $l=\sum\limits_{k=1}^j|y_{m_k}-y_{n_k}|$. Since for $m_k\geq
n_k+1$
\begin{equation}\label{(2.4)}
|\sum\limits_{i=n_k}^{m_k-1}\triangle y_i|
\leq\sum\limits_{i=n_k}^{m_k-1}|\triangle
y_i|\leq\sqrt{m_k-n_k-1}\left(\sum\limits_{i=n_k}^{m_k-1}|\triangle
y_i|^2\right)^{\frac{1}{2}},
\end{equation}
and since $(1-\cos y_i)\geq0$ and $y_i\notin
B_{\varepsilon}(\mathbf{\Theta})$ for  $i\in Z(n_k,m_k)$, we have
\begin{eqnarray}
J(y)&\geq&\frac{1}{2}\sum\limits_{k=1}^j\sum\limits_{i=n_k}^{m_k-1}|\triangle
y_i|^2+\sum\limits_{k=1}^j\sum\limits_{i=n_k}^{m_k-1}A(1-\cos
 y_i)\nonumber\\
&\geq&\sum\limits_{k=1}^j\left(\frac{|y_{m_k}-y_{n_k}|^2}{2r_k}+A\alpha_{\varepsilon}r_k\right)\nonumber\\
&\geq&\sum\limits_{k=1}^j\sqrt{{2A\alpha_{\varepsilon}}}|y_{m_k}-y_{n_k}|,\nonumber
\end{eqnarray}
where $r_k:=\begin{cases}
             m_k-n_k-1 &\text{ if } m_k>n_k+1\\
1  &\text{ if } m_k=n_k+1     \end{cases}$.
\end{proof}
\vs

Assume $0<\varepsilon<\gamma$. We will prove the existence of an
orbit minimizing function $J$ restricted to
$\Gamma_{\varepsilon}(\xi)$. For this purpose, we need the following
two lemmas.
\begin{lemma}\label{4}
Consider $\xi\in\mathbf{\Theta}\setminus\{0\}$ and assume that
$\{y^m\}_{m=1}^{\infty}\subset H$ is a minimizing sequence for
(\ref{(2.3)}) restricted to $\Gamma_{\varepsilon}(\xi)$, such that
for any $n\in{\mathbb N}, y^m\rightarrow y$  uniformly for  $i\in
Z[-n,n]$. If $y\in H$ and $J(y)<\infty$, then $y\in
\Gamma_{\varepsilon}(\xi)$.
\end{lemma}
\begin{proof}
By Remark \ref{R3}, there exist $\zeta$, $\eta\in\mathbf{\Theta}$
such that $y_{-\infty}=\zeta,\ y_{+\infty}=\eta$. By assumption
$y^m\rightarrow y$ uniformly for $i\in Z[-n,n]$ and
$y^m\in\Gamma_{\varepsilon}(\xi)$. \vs

\noi{\bf Claim 1:} $y_n\notin
B_{\varepsilon}(\mathbf{\Theta}\setminus\{0,\xi\})$ for all
$n\in{\mathbb
N}$. \\
Indeed, if there exist  $n_0\in \bn$ and $\theta\in
\mathbf{\Theta}\setminus\{0,\xi\}$ such that $y_{n_0}\in
B_{\varepsilon}(\theta)$, then
$\delta:=|y_{n_0}-\theta|<\varepsilon$. Since $y^m\rightarrow y$
 uniformly for $i\in Z[-n_0,n_0]$, we have for sufficiently large $m$ that
 $|y^m_{n_0}-y_{n_0}|<\varepsilon-\delta$ and
$|y^m_{n_0}-\theta|\leq|y^m_{n_0}-y_{n_0}|+|\theta-y_{n_0}|<\varepsilon$, which is
 a contradiction.
 \vs

\noi{\bf Claim 2:}  $y_{\pm\infty}\in\{0,\xi\}$.\\
If $y_{-\infty}=\zeta\in \mathbf{\Theta}\setminus\{0,\xi\}$, then
for any $0<\varepsilon_1\leq\varepsilon$, $\exists N_1\in{\mathbb
N}$ $\forall\ n\geq N_1$ $y_{-n}\in B_{\varepsilon_1/2}(\zeta)$.
Since $y^m\rightarrow y$ uniformly for $n\in Z[-N_1,N_1]$, there
exists $N_2\in{\mathbb N}$, $|y_{-N_1}^m-y_{-N_1}|<\varepsilon_1/2$
for $\forall\ m>N_2$.  Consequently, for those $\varepsilon_1$,
$N_1$, $N_2$ and $m>N_2$, we have
$|y_{-N_1}^m-\zeta|\leq|y_{-N_1}^m-y_{-N_1}|+|y_{-N_1}-\zeta|<\varepsilon_1\leq\varepsilon$.
Thus $y_{-N_1}^m\in B_{\varepsilon}(\zeta)$, which contradicts the
fact that $y^m\in \Gamma_{\varepsilon}(\xi)$. Thus
$\zeta\in\{0,\xi\}$. A similar argument can be applied to show
$\eta\in\{0,\xi\}$. \vs

\noi{\bf Claim  3:} $y_{-\infty}=0$.\\
 Since $y^m\in \Gamma_{\varepsilon}(\xi)$, for every $m\in{\mathbb N}$, there exists $n(m)\in{\mathbb Z}$ such
that $y_{n(m)+1}^m\notin B_{\varepsilon}(0)$ and  $y_n^m\in
B_{\varepsilon}(0)$ for all $n\leq n(m)$. For $y\in H$, put
$x_n(m):=y_{n-m}$ and   $x(m)=\{x_n(m)\}$. Then $J(x(m))=J(y)$.
Therefore, we can assume that $n(m)=0$ for all $m\in{\mathbb N}$.
Consequently $y_n^m\in B_{\varepsilon}(0)$ and $y_n\in
\overline{B}_{\varepsilon}(0)$, $\forall\ n\leq0$. Thus,
$\zeta\in \overline{B}_{\varepsilon}(0)\cap\{0,\xi\}=\{0\}$, i.e.
$\zeta=0$. \vs

\noi{\bf Claim 4:} $y_{+\infty}=\xi$.\\
Notice that $y_{+\infty}\in\{0,\xi\}$. Choose $\delta>0$ satisfying
$6\delta<\varepsilon$ and
$\frac{1}{2}(2\delta)^2+\delta^2<\sqrt{2A\alpha_{\delta}}\varepsilon/6$.
In order to show that such $\delta$ exists, put
$f(x)=\sqrt{2A\alpha_{x}}\varepsilon/6-\frac{1}{2}(2x)^2-x^2
=\varepsilon\sqrt{A}\sin\frac{x}{2}/3-3x^2$. Then,
$f^{\prime}(x)=\varepsilon\sqrt{A}\cos \frac{x}{2}/6-6x$, $f(0)=0$
and there exists $x_0\in(0,\pi/2)$ such that $f(x)>0$ for $0<x<x_0$,
which implies the existence of $\delta$ with the required
properties. Suppose, to the contrary, that $y_{+\infty}=0$, then
there exists $n_0\in{\mathbb N}$ such that $\forall n>n_0$
\;$y_{n_0}\notin B_{\delta}(0)$ and $y_n\in B_{\delta}(0)$. Since
$y^m\rightarrow y$ uniformly for $i\in Z[-n_0-1,n_0+1]$, there
exists a sufficiently large $m$, such that
$|y_{n_0+1}^m-y_{n_0+1}|<\delta$. Thus $y_{n_0+1}^m\in
B_{2\delta}(0)$. We need to consider the following two cases: \vs

\noi\textbf{Case 1:}
$y_{n_0}^m\notin B_{5\delta}(0)$. \\
Then
$|y_{n_0+1}^m-y_{n_0}^m|>3\delta$, and we have
\begin{eqnarray}
J(y^m)\geq9\delta^2/2+\sum\limits_{n=n_0+1}^{\infty}
[\frac{1}{2}|\triangle y_n^m|^2+A(1-\cos y_n^m)].\nonumber
\end{eqnarray}
Define
\begin{eqnarray}
x_n^m:= \left\{
\begin{array}{ll}
   0 ,    &   n\leq n_0 \\
  y_n^m,  &   n\geq n_0+1
\end{array}\right. \nonumber
\end{eqnarray}
Then $x^m=\{x^m_n\}\in\Gamma_{\varepsilon}(\xi)$ and
\begin{eqnarray}
J(x^m)&=&\sum\limits_{n=n_0}^{\infty}\left[\frac{1}{2}|\triangle
x^m_n|^2+A(1-\cos x^m_n)\right]\nonumber\\
&=&\frac{1}{2}|y_{n_0+1}^m|^2+\sum\limits_{n=n_0+1}^{\infty}\left[\frac{1}{2}|\triangle
y^m_n|^2+A(1-\cos y^m_n)\right]\nonumber\\
&\leq&
\frac{1}{2}|y_{n_0+1}^m|^2+J(y^m)-\frac{9\delta^2}{2}\nonumber\\
&<& J(y^m)-\frac{5\delta^2}{2},\nonumber
\end{eqnarray}
which leads to the following contradiction
\[c_{\varepsilon}(\xi)=\lim\limits_{m\rightarrow\infty}J(y^m)\geq
\lim\limits_{m\rightarrow\infty}J(x^m)+5\delta^2/2\geq
 c_{\varepsilon}(\xi)+5\delta^2/2.\]\\
\vs

\noi\textbf{Case 2:} $y_{n_0}^m\in B_{5\delta}(0)$.\\
\textbf{Subcase I:} $y_n^m\notin B_{\delta}(0)$ for all $1\leq
n\leq n_0$.\\
Then
\begin{eqnarray}
J(y^m)\geq\sqrt{2A\alpha_{\delta}}\varepsilon/6+
\sum\limits_{n=n_0+1}^{\infty}\left[\frac{1}{2}|\triangle
y^m_n|^2+A(1-\cos y^m_n)\right]\nonumber
\end{eqnarray}
Define
\begin{eqnarray}
z_n^m:= \left\{
\begin{array}{ll}
   0,     &   n\leq n_0 \\
  y_n^m,  &   n\geq n_0+1
\end{array}\right. \nonumber
\end{eqnarray}
Then $z^m=\{z^m_n\}\in\Gamma_{\varepsilon}(\xi)$ and
\begin{eqnarray}
J(z^m)&=&\sum\limits_{n=n_0}^{\infty}\left[\frac{1}{2}|\triangle
z^m_n|^2+A(1-\cos z^m_n)\right]\nonumber\\
&=&\frac{1}{2}|y_{n_0+1}^m|^2+\sum\limits_{n=n_0+1}^{\infty}\left[\frac{1}{2}|\triangle
y^m_n|^2+A(1-\cos y^m_n)\right]\nonumber\\
&\leq&
\frac{1}{2}|y_{n_0+1}^m|^2+J(y^m)-\sqrt{2A\alpha_{\delta}}\varepsilon/6\nonumber\\
&<& J(y^m)-\delta^2,\nonumber
\end{eqnarray}
which yields the following contradiction
\[c_{\varepsilon}(\xi)=\lim\limits_{m\rightarrow\infty}J(y^m)\geq
\lim\limits_{m\rightarrow\infty}J(x^m)+\delta^2\geq
c_{\varepsilon}(\xi)+\delta^2.\]
\\
\textbf{Subcase II:}  There exists a $n_1\in Z[1,n_0]$ such that
$y_{n_1}^m\in B_{\delta}(0), y_n^m\notin B_{\delta}(0),\ \forall\
n\in Z[1,n_1-1]$.\\ If $y_{n_1}^m\notin B_{5\delta}(0)$, by a
similar argument as in Case 1, we get a contradiction. On the other
hand, if $y_{n_1}^m\in B_{5\delta}(0)$, then by the argument used in
Subcase I of Case 2, we again obtain  a contradiction. \vs
Consequently, $y\in \Gamma_\ve(\xi)$, which completes the proof.
\end{proof}
\vs

\begin{lemma}\label{5}
For any
$\varepsilon\in(0,\gamma),~\xi\in\mathbf{\Theta}\setminus\{0\}$,
there exists $y^0:=y(\varepsilon,\xi)\in\Gamma_{\varepsilon}(\xi)$
such that $J(y(\varepsilon,\xi))=c_{\varepsilon}(\xi)$, i.e.
$y(\varepsilon,\xi)$ minimizes $J|_{\Gamma_{\varepsilon}(\xi)}$.
\end{lemma}
\begin{proof}
Let $\{y^m\}_{m=1}^{\infty}$ be a minimizing sequence for
(\ref{(2.3)}). There exists a positive number $M>0$ such that $M\geq
J(y^m)\geq\frac{1}{2}\sum\limits_{n=-\infty}^{\infty}|\triangle
y^m_n|^2$. We claim that $\{y_0^m\}_{m=1}^{\infty}$ is a bounded
sequence. Suppose to the contrary that for any $k\in{\mathbb N}$
there exists  $m_k\in{\mathbb N}$ such that $|y_0^{m_k}|\geq k$.
Thus $\lim\limits_{k\rightarrow\infty}|y^{m_k}_0|=\infty$, and there
exists $k_0\in{\mathbb N}$ such that $y_0^{m_k}\notin
B_{\varepsilon}(\xi)$ when $k\geq k_0$. Consider $y_1^{m_k}$. \vs

\textbf{Case I:} If $y_1^{m_k}\in \overline{B}_{\ve}(\xi)$, then
$J(y^m)\geq |y_0^{m_k}-\xi-\varepsilon|^2/2$. Let
$k\rightarrow\infty$, we have $J(y^{m_k})\rightarrow\infty$, which
contradicts the assumptions. \vs

\textbf{Case II:} Otherwise, $y_1^{m_k}\notin
\overline{B}_{\ve}(\xi)$. Denote $n_k:=\{n>0:~y_{n+1}^{m_k}\in
B_{\varepsilon}(\xi),y_l^{m_k}\notin B_{\varepsilon}(\xi), \forall\
l\in Z[0,n]\}$. Then we have
\begin{equation}\label{(2.5)}
J(y^{m_k})\geq\sqrt{2A\alpha_{\varepsilon}}|y_0^{m_k}-y_{n_k}^{m_k}|+\frac{1}{2}|y_{n_k+1}^{m_k}-y_{n_k}^{m_k}|^2\
\mbox{for all}\ k>k_0.
\end{equation}
Let $k\rightarrow\infty$ in (\ref{(2.5)}), then
$|y_0^{m_k}-y_{n_k+1}^{m_k}|\rightarrow\infty$. But
$|y_0^{m_k}-y_{n_k+1}^{m_k}|\rightarrow\infty$ if and only if
$|y_0^{m_k}-y_{n_k}^{m_k}|+|y_0^{m_k}-y_{n_k+1}^{m_k}|\rightarrow\infty$
which is equivalent to
$\sqrt{2A\alpha_{\varepsilon}}|y_0^{m_k}-y_{n_k}^{m_k}|+\frac{1}{2}|y_{n_k+1}^{m_k}-y_{n_k}^{m_k}|^2\rightarrow\infty$,
which contradicts again the assumptions. \vs

Consequently, $\{y_0^{m_k}\}$ is a bounded sequence and, by the
definition of the norm on $H$, $\{y^{m_k}\}$ is a bounded sequence
in $H$. Therefore, passing to a subsequence if necessary, there is
$y^0\in H$ such that $y^m$ weakly converges to $y^0$ in $H$. \vs

We claim $J(y^0)<\infty$. Indeed, consider  $-\infty<s<t<\infty$ and
define  for  $y\in H$
$$J(s,t,y)=\sum\limits_{n=s}^t\left[\frac{1}{2}|\triangle
y_n|^2+A(1-\cos y_n)\right].$$ The  weak convergence of the sequence
$\{y^m\}$ to $y^0$ in the Hilbert space $H$ implies  that
$y_n^m\rightarrow y_n^0$ for any $n\in {\mathbb Z}$. Then,
$\{y_n^m\}_{n=s}^t$ converges uniformly to $\{y_n^0\}_{n=s}^t$.
Clearly, $J(s,t,y)$ is lower continuous, so  it is also  lower
semi-continuous. Combining $M\geq J(y^m)\geq J(s,t,y^m)$ with the
lower semi-continuity of $J(s,t,y)$, we have
\begin{equation}\label{(2.6)}
J(s,t,y^0)\leq\liminf\limits_{m\rightarrow\infty}J(s,t,y^m)\leq
c_{\varepsilon}(\xi)=\liminf\limits_{m\rightarrow\infty}J(y^m)\leq
M.
\end{equation}
Since $y^0\in H$ and $s,t$ are arbitrary, (\ref{(2.6)}) implies
$J(y^0)\leq\inf\limits_{y\in\Gamma_{\varepsilon}(\xi)}J(y)$. Lemma \ref{4}
implies $y^0\in \Gamma_{\varepsilon}(\xi)$, and we have
$J(y^0)=c_{\varepsilon}(\xi)$.
\end{proof}
\vs

Put
\begin{equation}\label{c-eps}
c_{\varepsilon}:=\inf\limits_{\xi\in\mathbf{\Theta}\setminus\{0\}}c_{\varepsilon}(\xi).
\end{equation}
We will show that, fixed $\ve>0$, there are finite
$\zeta(\varepsilon)$'s such that
$\zeta(\varepsilon)\in\mathbf{\Theta}\setminus\{0\}$,
$c_{\varepsilon}(\zeta(\varepsilon))$ is a critical value of $J$
restricted on the set
$\bigcup_{\xi\in\mathbf{\Theta}}\Gamma_{\varepsilon}(\xi)$. \vs

\begin{lemma}\label{6}
The set $\Upsilon_\ve:=\{ \xi\in \Theta\setminus \{0\}:
c_\ve(\xi)=c_\ve\}$ is finite.
\end{lemma}
\begin{proof} Consider
$\xi\in\mathbf{\Theta}\setminus \{0\}$ and
$y\in\Gamma_{\varepsilon}(\xi)$. Then $y_{-\infty}=0$,
$y_{+\infty}=\xi$, $y_n\notin
B_{\varepsilon}(\mathbf{\Theta}\setminus\{0,\xi\})$. Put
$m_1:=\max\{n : y_{n}\in B_{\varepsilon}(0),\ y_m\notin
B_{\varepsilon}(0), \forall~ m> n\}$ and $m_2=\min\{n : y_n\in
B_{\varepsilon}(\xi), n\geq m_1\}$. If $m_2> m_1+2$, then by Lemma
\ref{2},
\begin{eqnarray}
J(y)\geq \sum\limits_{n=m_1}^{m_2-1}\frac{1}{2}|\triangle y_n|^2
      \geq
  \sqrt{2A\alpha_{\varepsilon}}|y_{m_2-1}-y_{m_1+1}|+\frac{1}{2}|\triangle
                y_{m_1}|^2+\frac{1}{2}|\triangle y_{m_2-1}|^2.
\nonumber
\end{eqnarray}
Notice that $\xi\rightarrow\infty$ if and only if
$|y_{m_2-1}-y_{m_1+1}|+|\triangle y_{m_1}|+|\triangle
y_{m_2-1}|\rightarrow \infty$ which is equivalent to
$\sqrt{2A\alpha_{\varepsilon}}|y_{m_2-1}-y_{m_1+1}|+\frac{1}{2}|\triangle
                y_{m_1}|^2+\frac{1}{2}|\triangle
 y_{m_2-1}|^2\rightarrow\infty$. Thus
$J(y)\rightarrow\infty$ as $\xi\rightarrow\infty$. In the case
$m_1\leq m_2\leq m_1+2$, by a  similar (but even simpler) argument,
we obtain the same result. Consider $\xi_0\in
\mathbf{\Theta}\setminus\{0\}$. Then we have
$c_{\varepsilon}(\xi_0)\geq c_{\varepsilon}$ and there exists
$M_1>0$ such that
$\inf_{y\in\Gamma_{\varepsilon}(\xi)}J(y)>c_{\varepsilon}(\xi_0)$
for all $|\xi|>M_1$. Consequently, there are only finitely many
$c_{\varepsilon}(\xi)$ which can be equal to $c_{\varepsilon}$.
\end{proof}
\vs

Fixed $\varepsilon>0$, Lemma \ref{6} implies $c_{\varepsilon}$ is
achieved  at some points $\zeta(\varepsilon)\in \Upsilon_\ve$. Now
by choosing a sequence of $\varepsilon_k\rightarrow0$, we claim that
there exists a subsequence $\{\epsilon_j\}_{j=1}^{\infty}$ such
that, for sufficiently large $j$ the points $\zeta(\varepsilon_j)\in
\Upsilon_{\ve_j}$ are independent of $j$, i.e. we have the
following:
\begin{lemma}\label{7}
Suppose that $\ve_k$ is a decreasing sequence of positive numbers
such that  $\ve_k\to 0$ as $k\to\infty$. Then  there exists a
subsequence $\{\varepsilon_j\}_{j=1}^{\infty}$ such that, for
sufficiently large $j$ the points $\zeta(\varepsilon_j)\in
\Upsilon_{\ve_j}$ are independent of $j$.
\end{lemma}
\begin{proof}
Consider $\{c(\varepsilon_k)\}_{k=1}^{\infty}$. For any
$y\in\Gamma_{\varepsilon_k}(\eta)$, we have $y_n\notin
B_{\varepsilon_k}(\mathbf{\Theta}\setminus\{0,\eta\})$ and also
$y_n\notin
B_{\varepsilon_{k+1}}(\mathbf{\Theta}\setminus\{0,\eta\})$ for all
$n\in {\mathbb Z}$. Thus $y\in\Gamma_{\varepsilon_{k+1}}(\eta)$ and
consequently
$\Gamma_{\varepsilon_{1}}(\eta)\subset\Gamma_{\varepsilon_{2}}(\eta)\subset
\dots\subset
\Gamma_{\varepsilon_{k}}(\eta)\subset\Gamma_{\varepsilon_{k+1}}(\eta)\subset\cdots$.
By definition of $c_{\varepsilon_k}(\eta)$, we have
\begin{equation}\label{(2.7)}
c_{\varepsilon_k}(\eta)=\inf\limits_{y\in
\Gamma_{\varepsilon_k}(\eta)}J(y)\geq\inf\limits_{y\in
\Gamma_{\varepsilon_{k+1}}(\eta)}J(y)=c_{\varepsilon_{k+1}}(\eta).
\end{equation}
Thus $\{c_{\varepsilon_k}\}_{k=1}^{\infty}$ is monotone
non-increasing bounded sequence. By a similar argument to the one
used in the proof of Lemma \ref{6}, the sequence
$\{\zeta(\varepsilon_k)\}_{k=1}^{\infty}$ is bounded. Consequently,
it contains a convergent  subsequence
$\{\zeta(\varepsilon_j)\}_{j=1}^{\infty}$. Since the set
$\mathbf{\Theta}$  consists of isolated points,
$\zeta(\varepsilon_j)$ is a constant sequence for $j$ sufficiently
large.
\end{proof}
\vs

Since for sufficiently large $j$ the points $\zeta(\varepsilon_j)$
are independent of $j$, denote by $\zeta=\zeta(\varepsilon_j)$. By
Lemma \ref{5}, there exists $y(\varepsilon_j, \zeta)\in
\Gamma_{\ve_j}(\zeta)$ such that
$c_{\varepsilon_j}=J(y(\varepsilon_j, \zeta))$.
\begin{theorem}\label{th1}
For $j$ sufficiently large, $y(\varepsilon_j,\zeta)$ is a
heteroclinic solution joining $0$ and $\zeta$.
\end{theorem}
\begin{proof}
Put $y(j):=y(\varepsilon_j,\zeta)$. By the definition of
$\Gamma_{\varepsilon}(\zeta)$ and $H$, it is sufficient to show that
for large $j$, $y_{n}(j)\notin \partial
B_{\varepsilon_j}(\mathbf{\Theta}\setminus\{0,\zeta\})$ for all
$n\in {\mathbb Z}$. If not, there would exist a sequence $
\eta_k\in\mathbf{\Theta}\setminus\{0,\zeta\}$ and $n_k\in{\mathbb
Z}$ such that
\begin{eqnarray}
y_{n_k}(k)\in\partial B_{\varepsilon_k}(\eta_k)\ \mbox{and}\
y_{n}(k)\notin\partial B_{\varepsilon_k}(\eta_k),\ \forall\
n<n_k.\nonumber
\end{eqnarray}
By similar argument used in the proof of Lemma \ref{6},  $\{\eta_k\}$ is bounded. Passing
to a subsequence, if necessary, $\eta_k$ must be a constant sequence, i.e.
 $\eta_k=:\eta$. We have the following two possibilities:\vs
\noi \textbf{Case 1:} There is an increasing sequence of integers
$k^{\prime}\rightarrow\infty$ such that $y_n(k^{\prime})\notin
\overline B_{\varepsilon_j}(\zeta)$ for all $n<n_{k^{\prime}}$, or\vs
\noi \textbf{Case 2:} For every $j\in{\mathbb N}$ there is a
$m_k<n_k$ such that $y_{m_k}(k)\in\partial
B_{\varepsilon_k}(\zeta)$. \vs

If Case 1 occurs, define
\begin{eqnarray}
x_n(k^{\prime})= \left\{
\begin{array}{ll}
  y_n(k^{\prime}),  &   n\leq n_k^{\prime} \\
  \eta ,           &   n\geq n_k^{\prime}+1
\end{array}\right. \nonumber
\end{eqnarray}
Then $y(k^{\prime})\in\Gamma_{\varepsilon_j}(\eta)$ and
\begin{align*}
J(y(k^{\prime}))-J(x(k^{\prime}))
&=\sum\limits_{n=n_{k^{\prime}}}^{\infty}\left[\frac{1}{2}|\triangle
            y_n(k^{\prime})|^2+A(1-\cos y_n(k^{\prime}))\right]\\
&\hskip.5cm -\frac{1}{2}|\triangle
 x_{n_{k^{\prime}}}(k^{\prime})|^2+A(1-\cos x_{n_{k^{\prime}}}(k^{\prime}))\\
&=\sum\limits_{n=n_{k^{\prime}}}^{\infty}\left[\frac{1}{2}|\triangle
            y_n(k^{\prime})|^2+A(1-\cos
y_n(k^{\prime}))\right]-\frac{1}{2}(\varepsilon_{k^{\prime}})^2-
A(1-\cos\varepsilon_{k^{\prime}}).
\end{align*}
If there exists a $n_0>n_{k^{\prime}}$ such that $y_{n_0}\notin
B_{\gamma}(\mathbf{\Theta})$, then
$J(y(k^{\prime}))-J(x(k^{\prime}))\geq
3A/2-(\varepsilon_{k^{\prime}})^2/2-A(1-\cos\varepsilon_{k^{\prime}})$.
Otherwise, there exist two adjacent points such that
 the distance of them is larger than $\gamma$. Then we have
$J(y(k^{\prime}))-J(x(k^{\prime}))\geq\sum\limits_{n=n_k}^{\infty}|\triangle
y_n(k^{\prime})|^2/2-(\varepsilon_{k^{\prime}})^2/2-A(1-\cos\varepsilon_{k^{\prime}})
 >2\pi^2/9-(\varepsilon_{k^{\prime}})^2/2-A(1-\cos\varepsilon_{k^{\prime}})$.
Define $\alpha:=\min\{3A/2-(\varepsilon_{k^{\prime}})^2/2-A(1-\cos
\varepsilon_{k^{\prime}}),
2\pi^2/9-(\varepsilon_{k^{\prime}})^2/2-A(1-\cos\varepsilon_{k^{\prime}})\}>0$.
We have $c_{\varepsilon_{k^{\prime}}}=J(y(k^{\prime}))\geq
J(x(k^{\prime}))+\alpha\geq c_{\varepsilon_{k^{\prime}}}+\alpha.$
This is a contradiction. \vs

If Case 2 occurs, define
\begin{eqnarray}
z_n(k):= \left\{
\begin{array}{ll}
  y_n(k),    &   n\leq m_k \\
  \zeta,      &   n\geq m_k+1
\end{array}\right. \nonumber
\end{eqnarray}
Then $z(k)\in\Gamma_{\varepsilon_j}(\zeta)$ and
\begin{eqnarray}
J(y(k))-J(z(k))
   &=&\sum\limits_{n=m_k}^{\infty}\left[\frac{1}{2}|\triangle
         y_n(k)|^2+A(1-\cos y_n(k)\right]-\frac{1}{2}|\triangle
 z_{m_k}(k)|^2-A(1-\cos z_{m_{k}}(k))\nonumber\\
   &=&\sum\limits_{n=m_k}^{\infty}\left[\frac{1}{2}|\triangle
         y_n(k)|^2+A(1-\cos
 y_n(k)\right]-\frac{1}{2}{\varepsilon_k}^2-A(1-\cos \varepsilon_k)\nonumber
\end{eqnarray}
By applying a similar argument as in the Case 1, we get again a
contradiction.
\end{proof}
\vs

As we can see on Figure \ref{fig:DE}, every heteroclinic solution
join two adjacent points of the set $\{2k\pi+\pi:~k\in\mathbb{Z}\}$,
or, after translation, heteroclinic solution join two adjacent
points of the set $\{2k\pi:~k\in\mathbb{Z}\}$. Denote by $\Upsilon$
the set of $\zeta\in\mathbf{\Theta}$ such that there exist a
heteroclinic solution joining $0$ to $\zeta$. The above observing
gives $\Upsilon=\{-2\pi,~2\pi\}$, which will be proved strictly
below. Since $A(1-\cos x)$ is $2\pi$-periodic, we have
$J(y+2\pi)=J(y)$. This implies that, for any integer $k>0$, if there
exists a heteroclinic orbit joining $-2k\pi$ and $0$, there must
exists a heteroclinic orbit joining $0$ and $2k\pi$. Thus we need
only to consider heteroclinic orbits joining $0$ to $2k\pi$.
\begin{lemma}\label{8}
$\Upsilon=\{-2\pi,~2\pi\}$.
\end{lemma}
\begin{proof}
Following the above argument, we just consider heteroclinic
solutions joining $0$ and $2k\pi$, where $k$ is a positive integer.
Suppose, to the contrary, Theorem \ref{th1} implies that there exist
$\zeta=2k\pi\in \Upsilon$ where $k>1$. Lemma \ref{5} guarantees
existence of $y$ which minimizes $J|_{\Gamma_{\varepsilon}(\zeta)}$.
Denote $n_1:=\min\{m:~y_n\in B_{\varepsilon}(\zeta),\forall\ n\geq
m\}$. We have the following two cases: \vs

\noi \textbf{Case 1:}
there exists  $n_0$ such that $y_{n_0}=y_{n_1-1}-2(k-1)\pi$.\\
Define
\begin{eqnarray}
x_n= \left\{
\begin{array}{ll}
  y_n,                                  &   n\leq n_0-1 \\
  y_{n+(n_1-n_0)}-2(k-1)\pi ,          &   n\geq n_0
\end{array}\right. \nonumber
\end{eqnarray}
Then $x\in\Gamma_{\varepsilon}(2\pi)$ and
\begin{eqnarray}
J(y)-J(x)=\sum\limits_{n=n_0}^{n_1-2}\left[\frac{1}{2}|\triangle
 y_n|^2+A(1-\cos y_n)\right]\nonumber
\end{eqnarray}
If $n_1-2=n_0, J(y)-J(x)\geq 4(k-1)^2\pi^2$. Otherwise, there exists
at least a suffix $n\in Z[n_0,n_1-2]$. If $n^{\prime}\in
Z[n_0,n_1-2]$ such that $y_{n^{\prime}}\notin
B_{\gamma}(\mathbf{\Theta})$, then we have
$J(y)-J(x)\geq\frac{3}{2}A$. Otherwise, there must be two adjacent
points such that the distant larger than $\gamma$. And then
$J(y)-J(x)>\frac{2\pi^2}{9}$. Define $\beta:=\min\{4(k-1)^2\pi^2,
3A/2, 2\pi^2/9\}$. All those situations contrary with
$c_{\varepsilon}(\zeta)=c_{\varepsilon}\geq c_{\varepsilon}(2\pi)+\beta\geq c_{\varepsilon}+\beta$. \\
\vs

\noi\textbf{Case 2:} If there is no $n_0$ such that
$y_{n_0}=y_{n_1-1}-2(k-1)\pi$, denote
$n_2:=\max\{n:~y_n<y_{n_1-1}-2(k-1)\pi\}$ and two situation maybe
meet: \vs

\noi\textbf{Subcase I:} If $n_2, n_1-1$ are two adjacent suffix.\\
Define
\begin{eqnarray}
x_n= \left\{
\begin{array}{ll}
  y_n                       &   n\leq n_2 \\
  y_{n}-2(k-1)\pi           &   n\geq n_2+1
\end{array}\right. \nonumber
\end{eqnarray}
Then $x\in\Gamma_{\varepsilon}(2\pi)$ and
\begin{eqnarray}
J(y)-J(x)= \frac{1}{2}|\triangle y_{n_2}|^2-\frac{1}{2}|\triangle
                 x_{n_2}|^2> 2(k-1)^2\pi^2, \nonumber
\end{eqnarray}
which implies the following contradiction
$$c_{\varepsilon}(\zeta)=c_{\varepsilon}\geq c_{\varepsilon}(2\pi)+\beta\geq c_{\varepsilon}+\beta.$$ \vs

\noi\textbf{Subcase II:} $n_2<n_1-2$. Then, we have
$y_{n_2}<y_{n_1-1}-2(k-1)\pi<y_{n_2+1}$.\\
Define
\begin{eqnarray}
x_n= \left\{
\begin{array}{ll}
  y_n                                 &   n\leq n_2 \\
  y_{n+(n_1-n_0)}-2(k-1)\pi           &   n\geq n_2+1
\end{array}\right. \nonumber
\end{eqnarray}
Then $x\in\Gamma_{\varepsilon}(2\pi)$ and
\begin{eqnarray}
J(y)-J(x)&=& \sum\limits_{n=n_2}^{n_1-2}[\frac{1}{2}|\triangle
             n|^2+A(1-\cos y_n)]-\frac{1}{2}|\triangle x_{n_2}|^2-A(1-\cos x_{n_2})\nonumber\\
        &=&  \sum\limits_{n=n_2+1}^{n_1-2}[\frac{1}{2}|\triangle
             y_n|^2+A(1-\cos y_n)]+\frac{1}{2}|\triangle y_{n_2}|^2
            -\frac{1}{2}|y_{n_1-1}-2(k-1)\pi-y_{n_1}|^2\nonumber
\end{eqnarray}
A similar argument as Case 1 of Theorem 1 induces a contradiction.

Consequently, we finish our proof.
\end{proof}
\begin{theorem}
For each $\xi\in\mathbf{\Theta}$, there exist at least two
heteroclinic orbits joining $\xi-2\pi$ to $\xi$ and at least two of
heteroclinic orbits joining $\xi$ to $\xi+2\pi$.
\end{theorem}
\begin{proof}
Without loss generality, we only need to check heteroclinic orbits
joining $0$ and $\zeta\in \Upsilon$. Lemma \ref{8} implies that only
$-2\pi$ and $2\pi$ belong to $\Upsilon$. If $\{y_n\}$ is a
heteroclinic orbit connecting $0$ and $2\pi$, then $\{y_{-n}\}$ is
also a heterclinic solution joining $2\pi$ to $0$. And
$\{y_n-2\pi\},\{y_{-n}-2\pi\}$ also two heteroclinic solutions
joining $-2\pi$ to $0$. The proof is complete.
\end{proof}
\vs

\section{Reasons for choosing Phase Plane of (\ref{(1.1)}) with $A=0.1$}
For simplicity, we paint phase plane of (\ref{(1.2)}) with $A=1$ in
section 1. We should paint phase plane of (\ref{(1.1)}) with $A=1$
to compare with that of (\ref{(1.2)}). However, phase plane of
(\ref{(1.1)}) with $A=1$ (figure \ref{fig:ADE}) is so different from
that of (\ref{(1.2)}). Non-periodic solutions move between the upper
and lower half plane of (\ref{(1.1)}). At first glance, the phase
plane of (\ref{(1.1)}) is different from that of (\ref{(1.2)}) in
essence. But it is not. Those phenomena appear because of
approximation error. Approximation error depends on amplitude. When
amplitude $A$ equal to $10$, we paint the phase plane of
(\ref{(1.1)}) as figure \ref{fig:LADE}. All periodic solutions and
non-periodic solutions become disordered. That is why we choose the
phase plane of (\ref{(1.1)}) with $A=0.1$. \vs

\end{document}